\documentclass{article}
\usepackage{fullpage}
\usepackage{graphics,psfrag,epsfig}
\usepackage{amsfonts,latexsym,eucal,amsmath,amsthm,amssymb}

\begin{document}

\newcommand{\rum}{\rule{0.5pt}{0pt}}
\newcommand{\rub}{\rule{1pt}{0pt}}
\newcommand{\rim}{\rule{0.3pt}{0pt}}
\newcommand{\numtimes}{\mbox{\raisebox{1.5pt}{${\scriptscriptstyle \times}$}}}
\newcommand{\optprog}[2]
{%
  \noindent\mbox{}\\[0cm]
  \noindent\fbox{%
  \begin{minipage}{0.955\linewidth}
    \mbox{}\\[-0.5cm]
    #1\\[#2]
  \end{minipage}
  }
  \noindent\mbox{}\\[-0.2cm]
}

\renewcommand{\refname}{References}

\twocolumn[%
\begin{center}
{\Large\bf The Collatz 3n+1 Conjecture is Unprovable \rule{0pt}{13pt}}\par
\bigskip
Craig Alan Feinstein \\
{\small\it 2712 Willow Glen Drive, Baltimore, Maryland 21209\rule{0pt}{13pt}}\\
\raisebox{-1pt}{\footnotesize E-mail: cafeinst@msn.com, BS"D}\par
\bigskip\smallskip
{\small\parbox{11cm}{%
\bigskip \noindent \textbf{Abstract:} In this paper, we show that any proof of
the Collatz $3n+1$ Conjecture must have an infinite number of lines; therefore, no
formal proof is possible.

\bigskip \noindent \textbf{Disclaimer:} This article was authored
by Craig Alan Feinstein in his private capacity. No official
support or endorsement by the U.S. Government is intended or
should be inferred.\rule[0pt]{0pt}{0pt}}}\bigskip
\end{center}]{%

\bigskip\noindent In 2005, the famous mathematician Freeman Dyson was asked, 
``What do you believe is true even though you cannot prove it?"
He answered:

``Since I am a mathematician, I give a precise answer to this question. 
Thanks to Kurt G\"{o}del, we know that there are true mathematical statements 
that cannot be proved. But I want a little more than this. I want a statement 
that is true, unprovable, and simple enough to be understood by people who are not 
mathematicians. Here it is.

``Numbers that are exact powers of two are 2, 4, 8, 16, 32, 64, 128 and so on. Numbers 
that are exact powers of five are 5, 25, 125, 625 and so on. Given any number such as 
131072 (which happens to be a power of two), the reverse of it is 270131, with the same 
digits taken in the opposite order. Now my statement is: it never happens that the reverse 
of a power of two is a power of five.

``The digits in a big power of two seem to occur in a random way without any regular
pattern. If it ever happened that the reverse of a power of two was a power of five,
this would be an unlikely accident, and the chance of it happening grows rapidly
smaller as the numbers grow bigger. If we assume that the digits occur at random,
then the chance of the accident happening for any power of two greater than a billion
is less than one in a billion. It is easy to check that it does not happen for powers
of two smaller than a billion. So the chance that it ever happens at all is less than
one in a billion. That is why I believe the statement is true.

``But the assumption that digits in a big power of two occur at random also implies that
the statement is unprovable. Any proof of the statement would have to be based on some
non-random property of the digits. The assumption of randomness means that the
statement is true just because the odds are in its favor. It cannot be proved because
there is no deep mathematical reason why it has to be true. (Note for experts: this
argument does not work if we use powers of three instead of powers of five. In that
case the statement is easy to prove because the reverse of a number divisible by
three is also divisible by three. Divisibility by three happens to be a non-random
property of the digits).

``It is easy to find other examples of statements that are likely to be true but
unprovable. The essential trick is to find an infinite sequence of events, each of
which might happen by accident, but with a small total probability for even one of
them happening. Then the statement that none of the events ever happens is probably
true but cannot be proved." \cite{b:Dy05}

In the spirit of Dyson's observation, we shall give an example of a 
statement that is likely to be true and then take things one step further by presenting a formal
proof that the statement is unprovable. Consider the following function:

\bigskip\noindent \textbf{Definition 1:} \textit{Let $T:\mathbb N
\rightarrow \mathbb N$ be a function such that $T(n)=\frac{3n+1}{2}$
if $n$ is odd and $T(n)=\frac{n}{2}$ if $n$ is even.}

\bigskip The Collatz $3n+1$ Conjecture states that for each $n \in \mathbb N$,
there exists a $k \in \mathbb N$ such that $T^{(k)}(n)=1$, where
$T^{(k)}(n)$ is the function $T$ iteratively applied $k$ times to
$n$ \cite{b:La85}. As of May 10, 2011, this conjecture has been verified for
all positive integers up to about $2^{60}$ \cite{b:Ro03}. Furthermore, one can give a
heuristic probabilistic argument \cite{b:Cr78} that since every
iterate of the function $T$ decreases on average by a
multiplicative factor of about $(\frac{3}{2})^{1/2}
(\frac{1}{2})^{1/2}=$ $(\frac{3}{4})^{1/2}$, all iterates will
eventually converge into the infinite cycle $\{1,2,1,2,...\}$,
assuming that each $T^{(k)}$ sufficiently mixes up $n$ as if each
$T^{(k)}(n)$ (mod 2) were drawn at random from the set $\{0,1\}$.

However, the Collatz $3n+1$ Conjecture has never been formally
proven. We shall prove that the Collatz $3n+1$ Conjecture
can, in fact, never be formally proven, even though there is a lot
of evidence for its truth. The underlying assumption in our argument is that
any proof of a theorem can be written in a computer text-file, which is composed of bits
(zeroes and ones). First, let us present a definition of ``random".

\bigskip\noindent \textbf{Definition 2:} \textit{We shall say that
vector ${\rm {\bf x}} \in \{0,1\}^k$ is random if ${\rm {\bf x}}$ cannot be specified 
in less than $k$ bits in a computer text-file \cite{b:Ch90}.}

\bigskip\noindent For example, the vector of one million concatenations
of the vector $(0,1)$ is not random, since we can specify it
in less than two million bits in a computer text-file by just writing, ``the vector of one
million concatenations of the vector $(0,1)$" in the text-file. However,
the vector of outcomes of one million
coin-tosses has a good chance of fitting our definition of
``random", since much of the time the most compact way
of specifying such a vector is to simply make a list of the
results of each coin-toss, in which one million bits are
necessary. We now prove three theorems.

\bigskip\noindent \textbf{Theorem 1:} \textit{For any
vector ${\rm {\bf x}} \in \{0,1\}^k$, there exists an $n \in \mathbb N$ such that
${\rm {\bf x}}=(n,T(n),...,T^{(k-1)}(n))$} (mod 2).

\bigskip\noindent \textit{Proof:} A proof of this can be found in
``The $3x+1$ problem and its generalizations" \cite{b:La85}.\qed

\bigskip\noindent\textbf{Theorem 2:} \textit{If $k,n \in \mathbb N$ and $T^{(k)}(n)=1$, then
in order to prove that $T^{(k)}(n)=1$, it is necessary to specify the values of
$(n,T(n),...,T^{(k-1)}(n))$} (mod 2) \textit{in the proof.}

\bigskip\noindent \textit{Proof:} Let the vector $(x_0(n),x_1(n),...,x_{k-1}(n))$ equal
$(n,T(n),...,T^{(k-1)}(n))$ (mod 2).
Then notice that the formula, $T^{(k)}(n)=\lambda_k(n)n+\rho_k(n)$ \cite{b:La85}, where 
$$
\lambda_k(n)\,=\frac{3^{x_0(n)+...+x_{k-1}(n)}}{2^k}
$$ and
$$
\rho_k(n)\,=\displaystyle\sum_{i=0}^{k-1}x_i(n)\frac{3^{x_{i+1}(n)+...+x_{k-1}(n)}}{2^{k-i}},
$$
is determined by the values of $(n,T(n),...,T^{(k-1)}(n))$ (mod 2) and there is a 
one-to-one correspondence between all of the possible formulas for $T^{(k)}(n)$ and 
all of the possible values of $(n,T(n),...,T^{(k-1)}(n))$ (mod 2); 
therefore, in order to prove that $T^{(k)}(n)=1$, it is necessary to specify the 
values of $(n,T(n),...,T^{(k-1)}(n))$ (mod 2) in the proof, since in order to prove 
that $T^{(k)}(n)=1$, it is necessary to specify the formula for $T^{(k)}(n)$ in the 
proof. \qed

\bigskip\noindent \textbf{Theorem 3:} \textit{It is impossible
to prove the Collatz $3n+1$ Conjecture.}

\bigskip\noindent \textit{Proof:} Suppose that there exists a proof of the
Collatz $3n+1$ Conjecture, and let $L$ be the number of bits in such
a proof. Now, let ${\rm {\bf x}} \in \{0,1\}^{L+1}$ be a random
vector, as defined above. (It is not difficult to prove that at
least half of all vectors in $\{0,1\}^{L+1}$ are random \cite{b:Ch90}.)
By Theorem 1, there exists an $n \in \mathbb N$ such that
${\rm {\bf x}}=(n,T(n),...,T^{(L)}(n))$ (mod 2) and
$T^{(L+1)}(n)=T^{(L)}(n)$ (mod 2). Then $T^{(L)}(n)>2$, so if
$T^{(k)}(n)=1$, then $k>L$. Hence, by Theorem 2 it is necessary to
specify the values of $(n,T(n),...,T^{(L)}(n))$ (mod 2) in order to
prove that there exists a $k \in \mathbb N$ such that
$T^{(k)}(n)=1$. But since $(n,T(n),...,T^{(L)}(n))$ (mod 2) is a
random vector, at least $L+1$ bits are necessary to specify
$(n,T(n),...,T^{(L)}(n))$ (mod 2), contradicting our assumption that
the proof contains only $L$ bits; therefore, a formal proof of the
Collatz $3n+1$ Conjecture cannot exist. \qed

}


\smallskip
\begin{thebibliography}{99}\small

\bibitem{b:Dy05} Dyson, F., ``What do you believe is true even though you cannot
prove it?", 
http://www.edge.org/q2005/q05\_9.html

\bibitem{b:La85} Lagarias, J.C., ``The $3x+1$ problem and its generalizations",
\textit{Amer. Math. Monthly 92} (1985) 3-23. [Reprinted in:
\textit{Conference on Organic Mathematics}, Canadian Math. Society
Conference Proceedings vol 20, 1997, pp. 305-331],
http://www.cecm.sfu.ca/organics/papers

\bibitem{b:Ro03} Roosendaal, E. (2011+), ``On the $3x + 1$ problem",
http://www.ericr.nl/wondrous/index.html

\bibitem{b:Cr78} Crandall, R. E., ``On the $3x+1$ problem",
\textit{Math. Comp.}, 32 (1978) 1281-1292.

\bibitem{b:Ch90} Chaitin, G.J., \textit{Algorithmic Information Theory},
revised third printing, Cambridge University Press, 1990.

\end{thebibliography}
\end{document}